\documentclass[a4paper]{article}

\usepackage[latin1]{inputenc}
\usepackage[T1]{fontenc}
\usepackage[francais,english]{babel}

\usepackage{amsmath,amssymb,amsfonts,amsthm}
\usepackage{amsfonts}
\usepackage{amssymb}
\usepackage{latexsym}
\usepackage{comment}

\def\R{\mathbb R}
\def\N{\mathbb N}

\def\et0{e^{tA}x_0}

\newtheorem{thm}{Theorem}
\newtheorem{Prop}{Proposition}

\author{
K. \textsc{Beauchard} 
\footnote{CMLS, Ecole Polytechnique, 91128 Palaiseau Cedex, France,
email: Karine.Beauchard@math.polytechnique.fr (corresponding author)},
P. \textsc{Cannarsa}
\footnote{Universit\`a di Roma Tor Vergata, via della Ricerca Scientifica 1, 00133, Roma, Italy, 
email: cannarsa@axp.mat.uniroma2.it},
\thanks{This research has been performed in the framework of the GDRE CONEDP.  The first author was partially supported by 
the ``Agence Nationale de la Recherche'' (ANR) Projet Blanc EMAQS number ANR-2011-BS01-017-01.}
}

\title{Inverse coefficient problem for Grushin-type parabolic operators}
\date{}

\begin{document}

\maketitle

\begin{abstract}
The approach to Lipschitz stability for uniformly parabolic equations introduced  by Imanuvilov and Yamamoto in 1998 based on Carleman estimates, 
seems hard to apply to the case of Grushin-type operators studied in this paper. Indeed, such estimates  are still missing for parabolic operators 
degenerating in the interior of the space domain.  Nevertheless, we are able to prove Lipschitz stability results for inverse coefficient problems 
for such operators, with locally distributed measurements in arbitrary space dimension. For this purpose, we follow a strategy that combines 
Fourier decomposition and Carleman inequalities for certain heat equations with nonsmooth coefficients (solved by the Fourier modes).
\end{abstract}


\section{Introduction}

\subsection{Model}

The relevance of the Heisenberg group to quantum mechanics has long been acknowledged. 
Indeed,  it was recognized by Weyl~\cite{Weyl}  that the Heisenberg algebra generated by the 
momentum and position operators comes from a Lie algebra representation associated with a 
corresponding group---namely the Heisenberg group (Weyl group in the traditional language of physicists).  
In such a group, the role played by the so-called Heisenberg laplacian is absolutely central, being analogous 
to the standard laplacian in Euclidean spaces, see \cite{MR983366}. On an even larger scale, deep connections 
have been pointed out between the properties of subriemannian operators, like the Heisenberg laplacian, 
and other topics of interest to current  mathematical research such as  isoperimetric problems and systems theory, see, for instance, \cite{capogna}. 

Another important example of sublaplacian is the Grushin operator which 
takes the form
\begin{equation}\label{eq:introG}
Gu = -(\partial^2_x u+x^2\partial_y^2u)
\end{equation}
on the plane.
As a matter of fact, the Heisenberg laplacian and the Grushin operator 
are deeply related: the former can be transformed into the latter, and the 
corresponding heat kernels are connected by an integral map, see \cite{MR2929725}.

This paper is a part of a general project we are pursuing, which consists of investigating the possibility of 
extending the known  controllability, observability, and Lipschitz stability properties of the heat equation,  
to degenerate parabolic problems. On all such topics, several results are available for parabolic operators 
which degenerate at the boundary of the space domain in low dimension, see, 
for instance, \cite{CTY1, Cannarsa-V-M-ADE, Cannarsa-V-M-SIAM, Ala-Can-Fra, Can-Fra-Roc_2, Cannarsa-V-M-CRAS}.

In two space dimensions, a fairly complete analysis of Grushin operator is presented in \cite{Grushin} as far as 
controllability and observability are concerned, and generalized to the multimenddimensional case in \cite{Grushin_IP}.
The inverse source problem is treated in \cite{Grushin_IP}.
To the best of our knowledge, there are no results on inverse coefficient problems for Grushin-type equations. 
The goal of this article is to prove a Lipschitz stability estimate for the inverse coefficient problem,
by adapting the techniques developed in \cite{Grushin_IP} for the inverse source problem.
\\

We consider Grushin-type equations of the form
\begin{equation} \label{Grushin}
\left\lbrace \begin{array}{ll}
\partial_t u - \Delta_x u - |x|^{2\gamma} b(x)  \Delta_y u = 0\,, & (t,x,y) \in (0,T)\times\Omega\,,\\
u(t,x,y)=0\,,                                                     & (t,x,y) \in (0,T)\times\partial \Omega\,,\\
u(0,x,y)=u^0(x,y)\,,                                              & (x,y) \in \Omega\, ,
\end{array}\right.
\end{equation}
where  
$T>0$,
$\Omega:=\Omega_1 \times \Omega_2$, 
$\Omega_1$ is a bounded open subset of $\mathbb{R}^{N_1}$, with $C^4$ boundary, such that $0 \in \Omega_1$,
$\Omega_2$ is a bounded open subset of $\mathbb{R}^{N_2}$, with $C^2$ boundary,
$N_1, N_2 \in \mathbb{N}^*:= \{1,2,3, ....\}$,
$b \in C^1(\overline{\Omega_1};(0,\infty))$, 
$\gamma \in (0,1]$  and $|.|$ is the Euclidean norm on $\mathbb{R}^{N_1}$.
\\

Specifically, we are interested in the inverse coefficient problem: 
is it possible to recover the coefficient $b(x)$ 
from the knowledge of an observation $\partial_t u |_{(T_0,T_1) \times \omega}$, 
where $\omega$ is a nonempty open subset of $\Omega$?

First, we recall well-posedness and regularity results for such equations.
To this aim, we introduce the space $H^1_\gamma(\Omega)$,
which is the closure of $C^\infty_0(\Omega)$ for the topology defined by the norm
$$\|f\|_{H^1_\gamma}:= \left( \int_{\Omega} \left( |\nabla_x f|^2 + |x|^{2\gamma} |\nabla_y f|^2 \right) dx dy \right)^{1/2},$$
and the Grushin operator $G_\gamma$ defined by
$$\begin{array}{ll}
&D(G_\gamma):=  
\{ f \in H^1_\gamma(\Omega); \exists c>0 \text{ such that }\\
&\left\vert 
\int_{\Omega} \left( \nabla_x f \cdot\nabla_x g 
+ |x|^{2\gamma} \nabla_y f\cdot
\nabla_y g \right) dx dy \right\vert \leqslant c \|g\|_{L^2(\Omega)}
\quad \text{for all $g \in H^1_\gamma (\Omega)$}\},
\end{array}
$$
$$
G_\gamma u := - \Delta_x u - |x|^{2\gamma} b(x) \Delta_y u.
$$

\begin{Prop} \label{Prop:WP}
Let $\gamma>0$.
For every $u_0 \in L^2(\Omega)$ and $g \in L^2((0,T)\times \Omega)$,
there exists a unique weak solution $u \in C^0([0,T];L^2(\Omega)) \cap L^2(0,T;H^1_\gamma(\Omega))$ of (\ref{Grushin}) such that
\begin{equation} \label{IC}
\begin{array}{ll}
u(0,x,y)=u_0(x,y) & (x,y) \in \Omega.
\end{array}
\end{equation}
Moreover, $u \in C^0((0,T];D(G_\gamma))$.
\end{Prop}

We refer to \cite{Grushin} for the proof with $N_1=N_2=1$; the general case 
can be treated similarly.

\subsection{Hypotheses and notations}

We introduce an open subset  $\Omega_1' \subset \subset \Omega_1$ such that $0 \notin \Omega_1'$
and $\delta>0$ such that
$$x\in \Omega_1' \quad \quad \Rightarrow \quad \quad |x| \geqslant \delta\,.$$
\\

The function $b$ is a priori assumed to satisfy
$$b \in \mathcal{M}:=\{ b \in C^0(\overline{\Omega_1},[m,M]) ; b \equiv 1 \text{ on } \Omega_1 \setminus \Omega_1'\}$$
for some positive constants $m, M$ with $0<m\leqslant 1 \leqslant M$.
Note that, in particular, $b \equiv 1$ on a neighborhood of $x=0$ and on a neighborhood of $\partial \Omega_1$.
\\

In order to introduce the hypotheses on the initial data $u^0$ of system (\ref{Grushin}), 
under which we prove Lipschitz stability estimate,
we first introduce several notations. Let $\mathcal{A}$ be the operator defined by
$$D(\mathcal{A}):=H^2 \cap H^1_0(\Omega_2), \quad \mathcal{A}\varphi := - \Delta_y \varphi,$$
let $(\mu_n)_{n \in \mathbb{N}^*}$ be the nondecreasing sequence of its eigenvalues
and $(\varphi_n)_{n \in \mathbb{N}^*}$ be the associated eigenvectors 
\begin{equation} \label{def:varphin}
\left\lbrace \begin{array}{ll}
- \Delta_y \varphi_n(y) = \mu_n \varphi_n(y)\,, & y \in \Omega_2\,,\\
\varphi_n(y)=0\,,                               & y \in \partial \Omega_2\,.
\end{array}\right.
\end{equation}
When $v=v(x,y) \in L^2(\Omega)$, then, we denote by $v_n=v_n(x)$ its Fourier components (with respect to variable $y$)
$$v_n(x):=\int_{\Omega_2} v(x,y) \varphi_n(y) dy, \forall n \in \mathbb{N}^*.$$

To prove the Lischitz stability estimate, 
the initial data $u^0$  of system (\ref{Grushin}) will be assumed to belong to the class
$$\begin{array}{ll}
\mathcal{D}_{N,K_1}:=\Big\{ &  u^0 \in D(G_\gamma^{s/2}) ;\, u_N^0 \geqslant 0 \text{ on } \Omega_1 \text{ and }  \\
                                & 
\underset{x \in \Omega_1'}{\sup} \Big( e^{t_1 \Delta_{\Omega_1'}} u_N^0\Big)(x)   \geqslant K_1 e^{\delta^{2\gamma} m T_1 \mu_N} \|u^0\|_{D(G_\gamma^{s/2})}  \Big\}
\end{array}$$
for some fixed time $t_1 \in (0,T_1)$,
some positive constants $s>N_1/2$, $K_1$ and some integer $N \in \mathbb{N}^*$.
Here, $e^{s \Delta_{\Omega_1'}}$ denotes the heat flow on $\Omega_1'$: 
for $\phi^0 \in L^2(\Omega_1')$, the function $\phi(s,x):=\left(e^{s \Delta_{\Omega_1'}}\phi^0\right)(x)$ is the solution of
$$\left\lbrace \begin{array}{ll}
\partial_s \phi(s,x) - \Delta_x \phi(s,x)=0\,, \quad & (s,x) \in (0,+\infty)\times \Omega_1'\,, \\
\phi(s,x)=0\,,                                       & (s,x) \in (0,+\infty) \times \partial \Omega_1'\,, \\
\phi(0,x)=\phi^0(x)                                  & x \in \Omega_1'\,.
\end{array}\right.$$
Note that
$$\Big\{  u^0 \in D(G_\gamma^{s/2}) ;\, u_N^0 \geqslant 0 \text{ on } \Omega_1 \Big\} = \cup_{j=1}^\infty \mathcal{D}_{N,1/j}.$$
In particular if $u^0 \in D(G_\gamma^{s/2})$ is $ \geqslant 0$ on $\Omega_1$ then 
$u^0_1 \geqslant 0$ thus $u^0 \in \mathcal{D}_{1,1/j}$ for $j$ large enough.
Thus, this class of functions is quite general.
\\

The letter $C$ denotes a constant that may change from line to line.

\subsection{Main results}

Our main result states the Lipschitz stability estimate, 
when the observation is done on a vertical strip $\omega:=\omega_1 \times \Omega_2$,
and for appropriate initial conditions.
When $\gamma \in(0,1)$, the Lipschitz stability estimate holds in any positive time.

\begin{thm} \label{thm:Main}
Let $\gamma \in (0,1)$, $0<T_1<T$, $\omega_1$ be a nonempty open subset of $\Omega_1$ and 
$\omega:=\omega_1 \times \Omega_2$ be a vertical strip.
There exists $$C=C(m,M,\Omega_1',K_1,t_1,T_1,\gamma,T,\omega_1)>0$$ such that,
for every $b, \tilde{b} \in \mathcal{M}$, $N \in \mathbb{N}$, $u^0 \in L^2(\Omega)$,
$\tilde{u}^0 \in \mathcal{D}_{N,K_1}$,
the associated solutions $u$ and $\tilde{u}$ of (\ref{Grushin}) satisfy
\begin{equation} \label{LSE}
\begin{array}{ll}
\displaystyle \int_{\Omega_1'} (b-\tilde{b})(x)^2  dx \leqslant 
& \displaystyle  \frac{C}{\|\tilde{u}^0\|_{D(G_\gamma^{s/2})}^2} \left[  \int_0^T \int_{\omega} |\partial_t (u-\tilde{u})(t,x,y)|^2 dx dy dt \right. \\
& \displaystyle + \left. \int_{\Omega_1' \times \Omega_2} |G_\gamma (u-\tilde{u})(T_1,x,y)|^2 dx dy \right]\,.
\end{array}
\end{equation}
\end{thm}

Note that the constant $C$ does not depend on $N$.
When $\gamma=1$, the Lipschitz stability estimate holds in time large enough, as stated below.

\begin{thm} \label{thm:Main_2}
We assume $\gamma=1$, $\omega_1$ be a nonempty open subset of $\Omega_1$ and 
$\omega:=\omega_1 \times \Omega_2$ be a vertical strip.
For $T_1>0$ large enough and all $T>T_1$
there exists $$C=C( m,M,\Omega_1',K_1,t_1,T_1,\gamma,T,\omega_1)>0$$ such that,
for every $b, \tilde{b} \in \mathcal{M}$, $N \in \mathbb{N}$, $u^0 \in L^2(\Omega)$,
$\tilde{u}^0 \in \mathcal{D}_{N,K_1}$,
the associated solutions $u$ and $\tilde{u}$ of (\ref{Grushin}) satisfy (\ref{LSE}).
\end{thm}

The assumption on the initial data ($u^0 \in  \mathcal{D}_{N,K_1,K_2}$ with $K_2 \geqslant \delta^{2\gamma} m T_1$)
is an important restriction, essentially related to technical difficulty. 
The validity of the Lipschitz stability estimate under more general assumptions is an interesting open problem
that will be investigated in future works.

\subsection{Structure of the article}

This article is organized as follows.
Section \ref{sec:Prel} is devoted to preliminary results concerning 
the well posedness of (\ref{Grushin}),
the Fourier decomposition of its solutions, 
the dissipation speed of the Fourier modes,
embeddings between spaces related to the Grushin operator,
and Harnack's inequality.
In Section \ref{sec:Main}, we prove our main results, namely, Theorems \ref{thm:Main} and \ref{thm:Main_2}.

\section{Preliminaries}
\label{sec:Prel}
\subsection{Well posedness}

\begin{Prop} \label{Prop:WP_ut}
Let $\gamma \in (0,1]$, $u^0 \in D(G_\gamma)$, $g \in H^1((0,T),L^2(\Omega))$, and
$u \in C^0([0,T],L^2(\Omega)) \cap L^2((0,T),H^1_\gamma(\Omega))$ be the solution of 
$$\left\lbrace \begin{array}{ll}
\partial_t u - \Delta_x u - |x|^{2\gamma} b(x)  \Delta_y u = g(t,x,y)\,, & (t,x,y) \in (0, \infty)\times\Omega\,,\\
u(t,x,y)=0\,,                                                            & (t,x,y) \in (0, \infty)\times\partial \Omega\,, \\
u(0,x,y)=u^0(x,y)\,,                                                     & (x,y) \in \Omega\,.
\end{array}\right.$$
Then the function $v:=\partial_t u$ belongs to $L^2((0,T),H^1_\gamma(\Omega))$ and solves
\begin{equation} 
\left\lbrace \begin{array}{ll}
\partial_t v - \Delta_x v - b(x) |x|^{2\gamma} \Delta_y v = \partial_t g(t,x,y)\,,    & (t,x,y) \in (0, \infty)\times\Omega\,,\\
v(t,x,y)=0\,,                                                                         & (t,x,y) \in (0, \infty)\times\partial \Omega\,,\\
v(0,x,y)=-G_\gamma u_0(x,y) + g(0,x,y)\,,                                             & (x,y) \in \Omega\,.
\end{array}\right.
\end{equation}
\end{Prop}

\subsection{Fourier decomposition}

\begin{Prop}
Let $u^0 \in L^2(\Omega)$, $g \in L^2((0,T)\times \Omega)$ and  $u$ be the solution of 
$$\left\lbrace \begin{array}{ll}
\partial_t u - \Delta_x u - |x|^{2\gamma} b(x)  \Delta_y u = g(t,x,y) & (t,x,y) \in (0, \infty)\times\Omega\,,\\
u(t,x,y)=0                                                            & (t,x,y) \in (0, \infty)\times\partial \Omega\,,\\
u(0,x,y)=u^0(x,y)                                                     & (x,y) \in \Omega
\end{array}\right.$$
For every $n \in \mathbb{N}^*$, the function
$$u_n(t,x):=\int_{\Omega_2} u(t,x,y) \varphi_n(y) dy$$
belongs to $C^0([0,T],L^2(\Omega))$ and is the unique weak solution of
\begin{equation} \label{Grushin_n}
\left\lbrace \begin{array}{ll}
\partial_t u_n - \Delta_x u_n + \mu_n |x|^{2\gamma} b(x) u_n = g_n(t,x)         & (t,x) \in (0,T)\times \Omega_1,\\
u_n(t,x)=0                                                                 & t \in (0,T)\times\partial\Omega_1,\\
u_n(0,x)=u_{n,0}(x)                                                        & x \in \Omega_1,
\end{array}\right.
\end{equation}
where
$$g_n(t,x):=\int_{\Omega_2} g(t,x,y) \varphi_n(y) dy
\quad \text{ and } \quad
u_{0,n}(x)=\int_{\Omega_2} u_0(x,y) \varphi_n(y) dy.$$
\end{Prop}
See \cite{Grushin} for the proof.

\subsection{Dissipation speed}

We introduce, for every $n \in \mathbb{N}^*, \gamma>0$, the operator $G_{n,\gamma}$ defined 
on $L^2(\Omega_1)$ by
\begin{equation} \label{def:An}
\begin{array}{ll}
D(G_{n,\gamma}):=H^2 \cap H^1_0(\Omega_1)\,, 
&
G_{n,\gamma} u := -\Delta_x u + \mu_n |x|^{2\gamma} b(x) u .
\end{array}
\end{equation}
The smallest eigenvalue of $G_{n,\gamma}$ is given by
$$
\displaystyle\lambda_{n,\gamma} = \min \left\{
\frac{ \int_{\Omega_1} \left[ |\nabla v(x)|^2 + \mu_n |x|^{2\gamma} b(x) v(x)^2 \right] dx }{\int_{\Omega_1} v(x)^2 dx } ;
v \in H^1_0(\Omega_1) ,\ v\neq 0\right\}\, .
$$
The asymptotic behavior (as $n \rightarrow + \infty$) of $\lambda_{n,\gamma}$,
which quantifies the dissipation speed of the solutions of (\ref{Grushin}),
is given by the following proposition (see \cite{Grushin_IP} for a proof).

\begin{Prop} \label{Prop:1st_eigenvalue}
For every $\gamma>0$, there exists constants $c_*, c^* >0$ such that
$$ c_* \mu_n^{\frac{1}{1+\gamma}} \leqslant \lambda_{n,\gamma} 
\leqslant c^* \mu_n^{\frac{1}{1+\gamma}}, \qquad \forall\, n \in \mathbb{N}^*\,.$$
\end{Prop}

\subsection{Continuous embeddings}
\label{subsec:embedding}

The goal of this section is to state and prove a continuous embedding used in the proof of the main result.

\begin{Prop}
For every $s>N_1/2$, we have $D(G_{N,\gamma}^{s/2}) \subset L^\infty(\Omega_1)$ with continuous embedding.
\end{Prop}

\textbf{Proof:}
We prove the conclusion just when $s$ is an even positive integer. In this case, setting $k=s/2$, it suffices to  show that
\begin{equation}\label{eq:imbedding}
k\in\N,\;\;k> \frac {N_1}4\qquad\Longrightarrow\qquad D(G_{N,\gamma}^{k}) \subset L^\infty(\Omega_1).
\end{equation}
Let $k=1$. We have that
 \begin{equation*}
u\in D(G_{N,\gamma})\Leftrightarrow u\in H^2 \cap H^1_0(\Omega_1)
\end{equation*}
Therefore, $u$ is continuous for $N_1=1$. Moreover,
\begin{equation*}
 u\in 
\begin{cases}
\hspace{.cm}
W^{1,p}(\Omega_1)\,,\;\forall p>1&
\text{if}\;\; N_1=2
\\
\hspace{.cm}
W^{1,2^{*}}(\Omega_1)
&
\text{if}\;\;N_1>2
\end{cases} 
\end{equation*}
where
\begin{equation*}
 \frac 1{2^*}=\frac{1}2-\frac 1{N_1}.
\end{equation*}
So, $u$ is H\"older continuos in $\overline{\Omega}_1$ by Sobolev's embedding provided that $2^*>N_1$, that is, $N_1<4$. We have thus checked \eqref{eq:imbedding} for $k=1$.

Now, suppose $N_1\geq 4$ (so that $2^*\le N_1$), let $k=2$,  and take $u\in D(G_{N,\gamma}^2) $. Set
$v:=G_{N,\gamma}u$, and observe that $u$ satisfies the boundary value problem
\begin{equation*}
\begin{cases}
\hspace{.cm}
 -\Delta_x u + \mu_N |x|^{2\gamma} b(x) u=v(x) &x\in \Omega_1
\\
u(x)=0
&x\in \partial\Omega_1.\hspace{.cm}
\end{cases}
\end{equation*}
Moreover, 
$ v\in 
W^{1,2^*}(\Omega_1,\R^{N_1})
$. Therefore,
\begin{equation*}
 v\in 
\begin{cases}
\hspace{.cm}
L^p(\Omega_1)\,,\;\forall p>1&
\text{if}\;\;2^*=N_1=4
\\
\hspace{.cm}
L^{2^{**}}(\Omega_1)
&
\text{if}\;\; N_1>4.
\end{cases} 
\end{equation*}
Thus, owing to the $L^p$-regularity of solutions to elliptic equations with H\"older continuous coefficients, 
\begin{equation}\label{eq:Sobolev}
 u\in 
\begin{cases}
\hspace{.cm}
W^{2,p}(\Omega_1)\,,\;\forall p>1&
\text{if}\;\;N_1=4
\\
\hspace{.cm}
W^{2,2^{**}}(\Omega_1)
&
\text{if}\;\; N_1>4.
\end{cases} 
\end{equation}
The above inclusions imply that $u$ is smooth right away if $2^{**}>N_1$, that is, 
$N_1<6$. By a refinement of the above argument one obtains the embedding  in \eqref{eq:imbedding}. Indeed, 
for $N_1\ge 6$, \eqref{eq:Sobolev} yields
\begin{equation*}
 u\in 
\begin{cases}
\hspace{.cm}
W^{1,p}(\Omega_1)\,,\;\forall p>1&
\text{if}\;\;2^{**}=N_1=6
\\
\hspace{.cm}
W^{1,2^{***}}(\Omega_1)
&
\text{if}\;\; N_1>6.
\end{cases} 
\end{equation*}
This gives that $u$ is H\"older continuous for $N_1<2^{***}$, that is, $N_1<8$.
We have thus checked \eqref{eq:imbedding} for $k=1,2$. The general result follows by  iteration.
$\Box$

\subsection{Harnack's inequality}

In this section, we recall the Harnack inequality for the heat equation (see \cite{evans})

\begin{Prop} \label{prop:Harnack}
Let $U$ be an open subset of $\Omega_1$, $T>0$, 
$V \subset \subset U$ connected, $0<t_1<T_1 <T$, $U_T:=(0,T) \times U$.
There exists $C_H>0$ such that
and for every $u \in C^2(U_T)$ with $u \geqslant 0$ on $U_T$ and
$$\partial_t u - \Delta u = 0 \text{ in } U_T$$
then
$$\underset{x \in V}{\inf} u(T_1,x) \geqslant C_H \underset{x \in V}{\sup} u(t_1,x).$$
\end{Prop}

\section{Proof of  Lipschitz stability }
\label{sec:Main}

In this section, first, we prove Theorem \ref{thm:Main}, then we explain how to adapt the reasoning to the proof of Theorem \ref{thm:Main_2}.
\\

The function $v(t,x,y):=(u-\tilde{u})(t,x,y)$ satisfies
$$\left\lbrace \begin{array}{ll}
\partial_t v_N - \Delta_x v_N + \mu_N |x|^{2\gamma} b(x) v_N = \mu_N |x|^{2\gamma} (b-\tilde{b})(x) \tilde{u}_N\,, & (t,x) \in (0,T)\times\Omega_1\,,\\
v_N(t,x)=0\,, & (t,x) \in (0,T)\times\partial \Omega_1\,, \\
v_N(0,x)=v_N^0(x)\,, & x \in \Omega_1\,.
\end{array}\right.$$

\textbf{Step 1: Use of Harnack inequality and assumption $\tilde{u}^0 \in \mathcal{D}_{N,K_1,K_2}$.}
We have
$$\begin{array}{ll}
\int_{\Omega_1'} (b-\tilde{b})(x)^2 dx
& \leqslant 
\frac{1}{\mu_N^2 \delta^{4\gamma} \underset{z \in \Omega_1'}{\inf} |\tilde{u}_N(T_1,z)|^2}
\int_{\Omega_1'} \left| \mu_N |x|^{2\gamma} (b-\tilde{b})(x) \tilde{u}_N(T_1,x) \right|^2  dx\,.
\end{array}$$
We recall that
$$\left\lbrace \begin{array}{ll}
\partial_t \tilde{u}_N - \Delta_x \tilde{u}_N + \mu_N |x|^{2\gamma} b(x) \tilde{u}_N = 0\,,        & (t,x) \in (0,+\infty)\times \Omega_1\,,\\
\tilde{u}_N(t,x)=0\,,                                                                              & t \in (0,+\infty)\times\partial\Omega_1\,,\\
\tilde{u}_N(0,x)=\tilde{u}_N^0(x)\,,                                                               & x \in \Omega_1\,.
\end{array}\right.$$
Let us introduce the solution $\nu_N(t,x)$ of
$$\left\lbrace \begin{array}{ll}
\partial_t \nu_N - \Delta_x \nu_N + \mu_N \delta^{2\gamma} m \nu_N = 0\,,  & (t,x) \in (0,T)\times\Omega_1' \,,\\
\nu_N(t,x)=0\,,                                                            & (t,x) \in (0,T)\times\partial \Omega_1'\,, \\
\nu_N(0,x)=\tilde{u}_N^0(x)\,,                                             & x \in \Omega_1' \,.
\end{array}\right.$$
Then
$$\begin{array}{c}
\mu_N |x|^{2\gamma} b(x) \geqslant \mu_N \delta^{2\gamma} m\,, \quad \forall x \in \Omega_1'\,,
\\
\tilde{u}_N(t,x) \geqslant 0 = \nu_N(t,x)\,, \quad \forall (t,x) \in (0,+\infty) \times \partial \Omega_1'\,,
\\
\tilde{u}_N(0,x) = \nu_N(0,x)\,, \quad \forall x \in \Omega_1'\,.
\end{array}$$
By the maximum principle, we deduce that
$$\tilde{u}_N(t,x) \geqslant \nu_N(t,x)\,, \quad \forall (t,x) \in (0,+\infty) \times \Omega_1'\,.$$
Note that 
$$\nu_N(t,x)= e^{- \mu_N \delta^{2\gamma} m t} \Big( e^{t \Delta_{\Omega_1'}} u_N^0 \Big)(x)\,.$$
Thus
$$\begin{array}{ll}
\underset{z \in \Omega_1'}{\inf} |\tilde{u}_N(T_1,z)|
& \geqslant \underset{z \in \Omega_1'}{\inf} |\nu_N(T_1,z)| \\
& \geqslant  e^{- \mu_N \delta^{2\gamma} m T_1} \underset{z \in \Omega_1'}{\inf} \Big| \Big( e^{T_1 \Delta_{\Omega_1'}} u_N^0 \Big)(z) \Big| \\
& \geqslant  e^{- \mu_N \delta^{2\gamma} m T_1} C_H \underset{z \in \Omega_1'}{\sup} \Big| \Big( e^{t_1 \Delta_{\Omega_1'}} u_N^0 \Big)(z) \Big| 
\text{ by Proposition \ref{prop:Harnack} } \\
& \geqslant 
C_H  K_1  \|u^0\|_{D(G_\gamma^{s/2})}
\end{array}$$
because $u^0 \in \mathcal{D}_{N,K_1}$. Therefore, we have
$$\begin{array}{ll}
\int_{\Omega_1'} (b-\tilde{b})(x)^2 dx
& \leqslant 
\frac{ 1 }{\mu_N^2 \delta^{4\gamma} C_H^2 K_1^2 \|u^0\|_{D(G_\gamma^{s/2})}^2 }
\int_{\Omega_1'} \left| \mu_N |x|^{2\gamma} (b-\tilde{b})(x) \tilde{u}_N(T_1,x) \right|^2  dx\,.
\end{array}$$
Thus
\begin{equation} \label{hyp}
\begin{array}{ll}
\int_{\Omega_1'} (b-\tilde{b})(x)^2 dx
& \leqslant 
\frac{C}{2\|u^0\|_{D(G_\gamma^{s/2})}^2 } \int_{\Omega_1'} \Big| \mu_N |x|^{2\gamma} (b-\tilde{b})(x) \tilde{u}_N(T_1,x) \Big|^2 dx
\\
& \leqslant \frac{C}{\|u^0\|_{D(G_\gamma^{s/2})}^2}  \int_{\Omega_1'}\Big( |\partial_t v_N(T_1,x)|^2 + |G_{N,\gamma} v_N(T_1,x)|^2 \Big) dx
\end{array}
\end{equation}
where 
$$C:= \frac{1}{\mu_1^2 \delta^{4\gamma} C_H^2 K_1^2}\,.$$
In order to dominate properly the first term of the right hand side, we revisit the proof of Proposition 6 of \cite{Grushin_IP}.
\\

\textbf{Step 2: Duhamel formula reads as}
$$\partial_t v_N(T_1)=e^{-G_{N,\gamma}(T_1-t)}\partial_t v_N(t)+\int_t^{T_1} e^{-G_{N,\gamma}(T_1-\tau)} g_N(\tau) d\tau, \quad \forall t \in (0,T_1)$$
where
$$g_N(\tau,x)=\mu_N |x|^{2\gamma} (b-\tilde{b})(x) \partial_t \tilde{u}_N(\tau,x).$$
Thus,
$$\begin{array}{ll}
\| \partial_t v_N(T_1) \|_{L^2(\Omega_1)} 
& \leqslant e^{-\lambda_{N,\gamma}(T_1-t)} \|\partial_t v_N(t)\|_{L^2(\Omega_1)} 
+ \int_t^{T_1} e^{-\lambda_{N,\gamma}(T_1-\tau)} \|g_N(\tau)\|_{L^2(\Omega_1)} d\tau\,.
\end{array}$$
Moreover,
$$\|g_N(\tau)\|_{L^2(\Omega_1)}
\leqslant C \mu_N  \| b-\tilde{b} \|_{L^2(\Omega_1')} \|\partial_t \tilde{u}_N(\tau)\|_{L^\infty(\Omega_1)}\,.$$
By the continuous embedding proved in Section \ref{subsec:embedding}, we have 
$$\begin{array}{ll}
\|\partial_t \tilde{u}_N(\tau)\|_{L^\infty(\Omega_1)}
& \leqslant
C \|\partial_t \tilde{u}_N(\tau)\|_{D(G_{N,\gamma}^{s/2})} \\
& \leqslant
C \| \tilde{u}_N^0 \|_{D(G_{N,\gamma}^{s/2})} e^{-\lambda_{N,\gamma} \tau} \\
& \leqslant
C \| \tilde{u}^0 \|_{D(G_{\gamma}^{s/2})} e^{-\lambda_{N,\gamma} \tau}.
\end{array}$$
Therefore,
\begin{equation} \label{interm2}
\| g_N(\tau) \|_{L^2(\Omega_1)}  \leqslant
 C \mu_N \| b-\tilde{b}\|_{L^2(\Omega_1')}\| \tilde{u}^0 \|_{D(G_{\gamma}^{s/2})} e^{-\lambda_{N,\gamma} \tau}
\end{equation}
and
$$\begin{array}{ll}
\| \partial_t v_N(T_1) \|_{L^2(\Omega_1)} \leqslant 
& e^{-\lambda_{N,\gamma}(T_1-t)} \|\partial_t v_N(t)\|_{L^2(\Omega_1)} \\
& + C (T_1-t) \| \tilde{u}^0 \|_{D(G_{\gamma}^{s/2})} \mu_N  e^{-\lambda_{N,\gamma}T_1} \| b-\tilde{b} \|_{L^2(\Omega_1')}\,.
\end{array}$$
Taking the square, we get
$$\begin{array}{ll}
\int_{\Omega_1} | \partial_t v_N(T_1,x) |^2 dx \leqslant 
& 2 e^{-2 \lambda_{N,\gamma}(T_1-t)} \int_{\Omega_1} | \partial_t v_N(t,x) |^2 dx  \\
& + 2 C^2 (T_1-t)^2 \| \tilde{u}^0 \|_{D(G_{\gamma}^{s/2})}^2 \mu_N^2  e^{-2 \lambda_{N,\gamma}T_1} \int_{\Omega_1'} (b-\tilde{b})(x)^2 dx.
\end{array}$$
Integrating over $t \in (T_1/3,2T_1/3)$, we obtain
\begin{equation} \label{Duham}
\begin{array}{ll}
\int_{\Omega_1} | \partial_t v_N(T_1) |^2  \leqslant 
& \frac{6}{T_1} e^{-2 \lambda_{N,\gamma}T_1/3} \int_{T_1/3} ^{2T_1/3} \int_{\Omega_1} | \partial_t v_N(t,x)) |^2 dx dt \\ 
& + 2 C^2 \left( \frac{2T_1}{3} \right)^2 \| \tilde{u}^0 \|_{D(G_{\gamma}^{s/2})}^2 \mu_N^2  e^{-2 \lambda_{N,\gamma}T_1} \int_{\Omega_1'} |x|^{4\gamma} (b-\tilde{b})(x)^2 dx.
\end{array}
\end{equation}

\textbf{Step 3: We apply Carleman estimate.} 
Working exactly as in the step 2 of the proof of Proposition 6 of \cite{Grushin_IP}, we get, for $N$ large enough
$$\begin{array}{ll}
\int_{T_1/3} ^{2T_1/3} \int_{\Omega_1} | \partial_t v_N(t,x) |^2 dx dt
\leqslant C e^{C \mu_N^{p(\gamma)}} \Big( 
& 
 \int_0^{T_1} \int_{\omega_1} |\partial_t v_N(t,x)|^2 dx dt \\
& +\int_0^{T_1} \int_{\Omega_1} |g_N(t,x)|^2 dx dt \Big)
\end{array}$$
where $C=C(T_1)>0$ and $p(\gamma)=1/2$ if $\gamma \in [1/2,1]$ and $p(\gamma)=2/3$ if $\gamma \in (0,1/2)$.
Moreover, estimate (\ref{interm2}) justifies that
$$\int_0^{T_1} \int_{\Omega_1} |g_N(t,x)|^2 dx dt
\leqslant 
\frac{C^2}{2\lambda_{N,\gamma}} \mu_N^2  \| \tilde{u}^0 \|_{D(G_{\gamma}^{3/2})}^2 \int_{\Omega_1'} (b-\tilde{b})(x)^2 dx.$$
Thus,
\begin{equation} \label{Carl}
\begin{array}{ll}
\int_{T_1/3} ^{2T_1/3} \int_{\Omega_1} | \partial_t v_N(t,x) |^2 dx dt \leqslant
& C e^{C \mu_N^{p(\gamma)}}   \int_0^{T_1} \int_{\omega_1} |\partial_t v_N(t,x)|^2 dx dt \\
& + \frac{C}{\lambda_{N,\gamma}} \mu_N^2 e^{C \mu_N^{p(\gamma)}}   \| \tilde{u}^0 \|_{D(G_{\gamma}^{s/2})}^2 \int_{\Omega_1'} (b-\tilde{b})(x)^2 dx.
\end{array}
\end{equation}

\textbf{Step 4: By combining (\ref{hyp}), (\ref{Duham}) and (\ref{Carl}), we obtain}
$$\begin{array}{ll}
 &  \int_{\Omega_1'} (b-\tilde{b})(x)^2 dx \\
\leqslant & \frac{C}{\|u^0\|_{D(G_\gamma^{s/2})}^2} \left(    \int_{\Omega_1'} |G_{N,\gamma} v_N(T_1,x)|^2 dx 
 +  e^{C \mu_N^{p(\gamma)}-2 \lambda_{N,\gamma}T_1/3}   \int_0^{T_1} \int_{\omega_1} |\partial_t v_N(t,x)|^2 dx dt \right)  \\
& +  \left(  
\frac{1}{\lambda_{N,\gamma}}  \mu_N^2 e^{C \mu_N^{p(\gamma)} -2 \lambda_{N,\gamma}T_1/3}  
+  \mu_N^2  e^{-2 \lambda_{N,\gamma}T_1}
\right) \int_{\Omega_1'} (b-\tilde{b})(x)^2 dx
\end{array}$$
 
For $\gamma \in (0,1)$, $0<T_1<T$ arbitrary and $N$ large enough, 
the source term in the right hand side may be absorbed by the left hand side 
(because $\mu_N^{p(\gamma)} << \lambda_{N,\gamma}$)
and we get a constant $\mathcal{C}>0$ (independent of $N$) such that
$$\begin{array}{ll}
          & \int_{\Omega_1'} (b-\tilde{b})(x)^2 |x|^{4\gamma} dx  \\
\leqslant & \frac{\mathcal{C}}{\|u^0\|_{D(G_\gamma^{s/2})}^2}  \left(\int_0^{T_1} \int_{\omega_1} |\partial_t v_N(t,x)|^2 dx dt + \int_{\Omega_1'} |G_{N,\gamma} v_N(T_1,x)|^2 dx \right) \\
\leqslant & \frac{\mathcal{C}}{\|u^0\|_{D(G_\gamma^{s/2})}^2}  \left(\int_0^{T_1} \int_{\omega} |\partial_t v(t,x,y)|^2 dx dy dt + \int_{\Omega} |G_{\gamma} v(T_1,x,y)|^2 dx dy \right) 
\end{array}$$

When $\gamma=1$, then $\lambda_{N,\gamma}$ behaves asymptotically like $C \mu_N^{p(\gamma)}$,
thus the time $T_1$ needs to be taken large enough for the same conclusion to hold.
\\

The previous arguments treat the high frequencies (i.e., $N \geqslant N_*$ for some $N_*$). 
For low frequencies (i.e., $N < N_*$), the Lipschitz stability estimate for the inverse source problem in the uniformly parabolic case (see \cite{YZ}) yields the conclusion.

\begin{center}
\textbf{Acknowledgements} 
\end{center}

The authors are grateful to Masahiro Yamamoto and Philippe Gravejat for fruitful discussions.

\bibliography{biblio}
\bibliographystyle{plain}

\end{document}